\documentclass{article}
\usepackage{amssymb,amsmath,epsfig,graphics}

\newtheorem{corollary}{Corollary}

\newtheorem{lemma}{Lemma}
\newtheorem{proposition}{Proposition}
\newtheorem{definition}{Definition}

\newtheorem{theorem}{Theorem}

\newcommand{\proof}{\textbf{Proof:  }}
\newcommand{\CC}{{\mathbb C }}
\newcommand{\FF}{{\mathbb F }}
\newcommand{\NN}{{\mathbb N }}
\newcommand{\PP}{{\mathbb P }}
\newcommand{\QQ}{{\mathbb Q }}

\newcommand{\ZZ}{{\mathbb Z }}
\newcommand{\RR}{{\mathbb R }}
\newcommand{\UU}{{\mathbb U }}

\newcommand{\WW}{{\mathbb W }}

\newcommand{\UL}{{\mathbb U }^\gL}

\newcommand{\gA}{\mathfrak{A}}

\newcommand{\cD}{\mathcal{D}}

\newcommand{\cI}{\mathcal{I}}

\newcommand{\sB}{\mathsf{B}}
\newcommand{\sD}{\mathsf{D}}
\newcommand{\osD}{\overline{\sD}}
\newcommand{\sH}{\mathsf{H}}
\newcommand{\sQ}{\mathsf{Q}}

\newcommand{\sF}{\mathsf{F}}
\newcommand{\mo}{\mathsf{m}}

\newcommand{\mul}{\mathrm{mult}}

\newcommand{\gG}{\Gamma}
\newcommand{\gD}{\Delta}
\newcommand{\gL}{\Lambda}
\newcommand{\vL}{\Lambda^\vee}
\newcommand{\gl}{\lambda}
\newcommand{\ep}{\epsilon}
\newcommand{\ga}{\alpha}

\newcommand{\ru}{\mu}
\newcommand{\rul}{\mu^\gL}

\newcommand{\qed}{{\hfill$\blacksquare$}}

\newcommand{\mf}[1]{\mathsf{#1}}
\newcommand{\wh}[1]{\widehat{#1}}

\newcommand{\FD}{\wh{\cD}}

\newcommand{\modquot}[2]{\mbox{\raisebox{.1em}{$#1$}\hspace{-2mm}
{ / }\hspace{-2mm} \raisebox{-.1em}{$#2$}}}

\begin{document}

\title{Motives from Diffraction.}

\author{Jan Stienstra\footnote{e-mail: stien{@}math.uu.nl}\\
\small{Mathematisch Instituut, Universiteit Utrecht, the Netherlands}\normalsize}

\date{}

\maketitle

\hfill\textit{Dedicated to Jaap Murre and Spencer Bloch}

\begin{abstract}
We look at geometrical and arithmetical patterns created from 
a finite subset of $\ZZ^n$ by diffracting waves and bipartite graphs.
We hope that this can make a link between Motives and 
the Melting Crystals/Dimer models in String Theory.
\end{abstract}

\section*{Introduction.}

\emph{
Why is it that, occasionally, mathematicians studying \emph{Motives} and
physicists searching for a \emph{Theory of Everything} seem to be looking at the same examples, just from different angles? Should the Theory of Everything include properties of Numbers? 
Does Physics yield realizations of Motives which have not been considered before in the cohomological set-up of motivic theory?}

Calabi-Yau varieties of dimensions $1$ and $2$, being elliptic curves
and K3-surfaces, have a long and rich history in number theory and geometry.
Calabi-Yau varieties of dimension $3$ have played an important role in many developments in String Theory. The discovery of \emph{Mirror Symmetry} attracted the attention of physicists and mathematicians to Calabi-Yau's near the \emph{large complex structure limit} \cite{M, Yau}. Some analogies between String Theory and Arithmetic Algebraic Geometry near this limit were discussed in \cite{S2, S3, S1}. Recently new models appeared, called \emph{Melting Crystals} and \emph{Dimers} \cite{ORV,KOS}, which led to interesting new insights in String Theory, without going near the large complex structure limit.
The present paper is an attempt to find motivic aspects of these new models.
We look at geometrical and enumerative patterns associated with 
a finite subset $\gA$ of $\ZZ^n$. The geometry comes from waves diffracting on $\gA$ and from a periodic weighted bipartite graph generated by $\gA$. The latter is related to the dimers
(although here we can not say more about this relation).
Since the tori involved in these models are naturally dual to each other there seems to be some sort of
mirror symmetry between the diffraction and the graph pictures. The enumerative patterns count lattice points on the diffraction pattern, points on varieties over finite fields
and paths on the graph. They are expressed through a sequence
of polynomials $\sB_{N}(z)$ with coefficients in $\ZZ$ and via a
limit for $z\in\CC$:
\begin{equation} \label{eq:radius}
\sQ(z)=\lim_{N\to\infty} |\sB_{N}(z)|^{-N^{-n}}\,.
\end{equation}
Limit formulas like (\ref{eq:radius})
appear frequently and in very diverse contexts in the literature, e.g.
for 
\emph{entropy in algebraic dynamical systems} in \cite{EW} Theorem 4.9, 
for 
\emph{partition function per fundamental domain in dimer models}
in \cite{KOS} Theorem 3.5, 
for \emph{integrated density of states} in \cite{GKT} p.206.
Moreover, $\sQ(z)$ appears as 
\emph{Mahler measure} in \cite{Bo}, as the exponential of a 
\emph{period in Deligne cohomology} in \cite{De, RV},
and in \emph{instanton counts} in \cite{S1}; see the remark at the end of Section \ref{section:geometry}.

With some additional restrictions $\gA$ provides the toric data for
a family of Calabi-Yau varieties and various well-known results about
Calabi-Yau varieties near the large complex structure limit can be derived from the Taylor series expansion of $\log\sQ(z)$ near $z=\infty$; see the Remark at the end of Section \ref{section:moments}.
In the present paper we are not so much interested in the large complex structure limit. Instead we focus on the polynomials $\sB_{N}(z)$ and the
limit formula (\ref{eq:radius}). 
\emph{This does not require conditions of `Calabi-Yau type'. }
\begin{figure}
\begin{center}
\setlength\epsfxsize{10cm}
\epsfbox{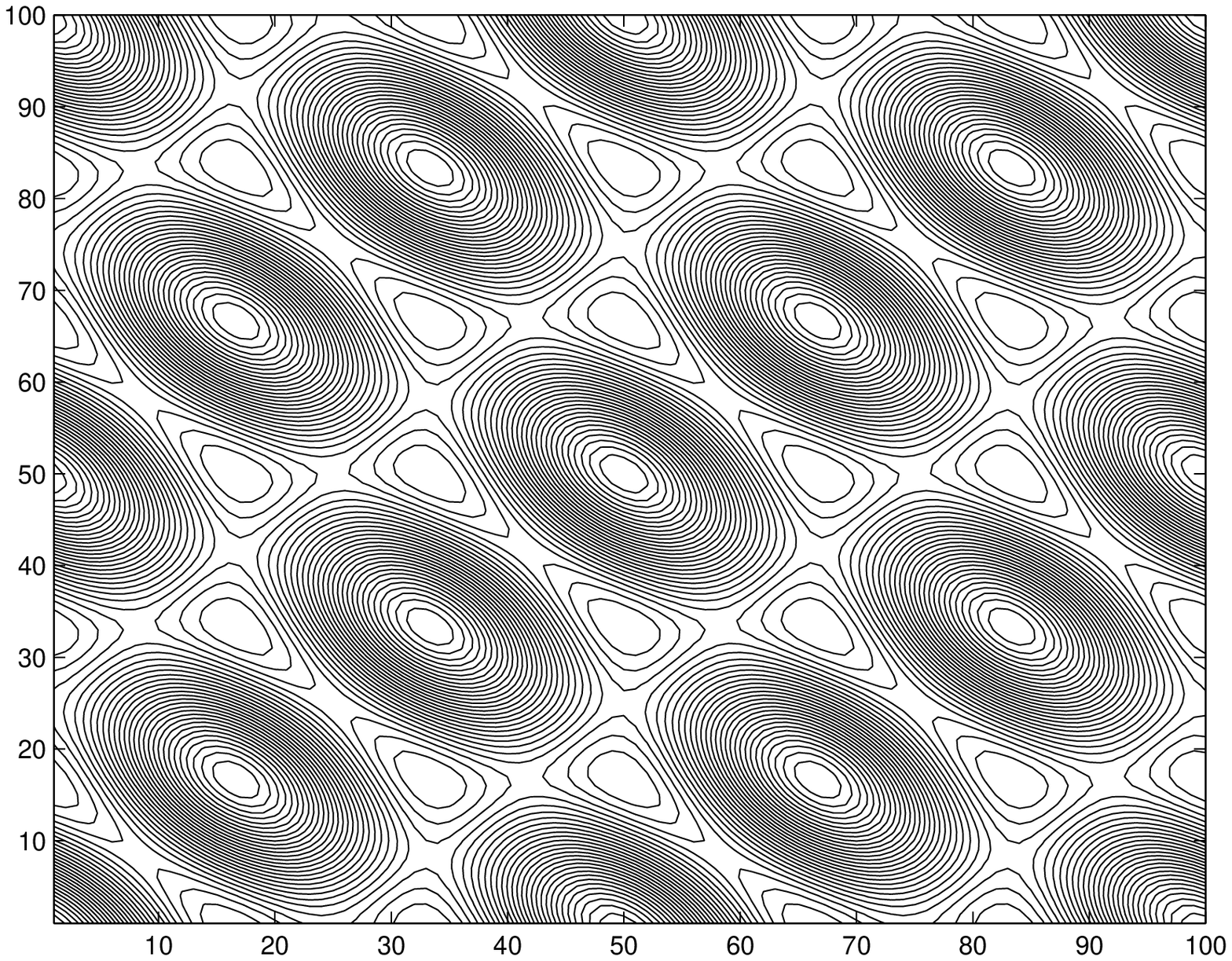}
\end{center}
\caption{\label{fig:diffraction pattern}
\textit{Diffraction pattern for
$\gA=\{(1,0),\,(0,1),\,(-1,-1)\}\subset\ZZ^2$, all 
$c_{\mf{a}}=1$.}}
\end{figure}

When waves are diffracted at some finite set $\gA$ of points in a plane, 
the diffraction pattern observed in a plane at large distance is, according to the Frauenhofer model, the absolute value squared of the Fourier transform of $\gA$.
There is no mathematical reason to restrict this model to dimension $2$. Also the points may have weights
$\geq 1$. So, we take a finite subset $\gA$ of $\ZZ^n$ and positive integers $c_\mf{a}\;(\mf{a}\in\gA)$. 
These data can be summarized as a distribution 
$
\cD=\sum_{\mf{a}\in\gA}c_\mf{a}\delta_\mf{a}\,,
$
where $\delta_\mf{a}$ denotes the Dirac delta distribution, evaluating test functions at the point $\mf{a}$.
The Fourier transform of $\cD$ is the function 
$
\FD(\mf{t})=\sum_{\mf{a}\in\gA}c_\mf{a} e^{-2\pi i\langle \mf{t},\mf{a}\rangle} 
$
on $\RR^n$;
here $\langle,\rangle$ is the standard inner product on $\RR^n$.
The diffraction pattern consists of the level sets of the function
$$
|\FD(\mf{t})|^2=\sum_{\mf{a},\mf{b}\in\gA}c_\mf{a}c_\mf{b} 
e^{2\pi i\langle \mf{t},\mf{a}-\mf{b}\rangle}
=\sum_{\mf{a},\mf{b}\in\gA}c_\mf{a}c_\mf{b} 
\cos(2\pi\langle \mf{t},\mf{a}-\mf{b}\rangle)\,.
$$
This function is periodic with period lattice $\vL$  
dual to the lattice $\gL$ spanned over $\ZZ$ by the differences
$\mf{a}-\mf{b}$ with $\mf{a},\mf{b}\in\gA$. \emph{Throughout this note we assume that $\gL$ has rank $n$.}
Looking at the intersections of the diffraction pattern with the lattices 
$\frac{1}{N}\vL$ we introduce the enumerative data
\begin{equation}\label{eq:multiplicities}
\mul_N(r):=\sharp\{\mf{t}\in\modquot{\textstyle{\frac{1}{N}}\vL}{\vL}\;|
\;|\FD(\mf{t})|^2\,=\,r\;\}
\qquad \textrm{for}\quad N\in\NN,\,r\in\RR\,,
\end{equation}
and use these to define polynomials $\sB_N(z)$ as follows:
\begin{definition}\label{def:BN}
\begin{equation}\label{eq:BN diff}
\sB_N(z):=\prod_{r\in\RR}(z-r)^{\mul_N(r)}\,.
\end{equation}
\end{definition}
One could also introduce the generating function
$
\sF (z,T):=\sum_{N\in\NN}\sB_N(z)\,T^{N^n}
$,
but except for the classical number theory of the case $n=1$ (see Section \ref{subsection:n=1}), and the observation that Formula (\ref{eq:radius}) gives $\sQ(z)$ as
the radius of convergence of $\sF (z,T)$ as a complex power series in $T$,
we do not yet have appealing results about $\sF (z,T)$.

For the graph model we start from the same data: the finite set $\gA\subset\ZZ^n$, the weights $c_\mf{a}\;(\mf{a}\in\gA)$ 
and the lattice $\gL$ spanned by the differences $\mf{a}-\mf{b}$ with $\mf{a},\,\mf{b}\in\gA$. We must now assume that
$$
\gA\cap\gL\,=\,\emptyset.
$$

\begin{figure}
\begin{picture}(300,160)(-90,-20)
\thinlines
\multiput(0,0)(0,75){2}{
\multiput(0,0)(-25,25){2}{
\multiput(0,0)(75,0){3}{
\begin{picture}(50,50)(0,0)
\put(0,0){\line(1,0){23}}
\put(0,0){\line(0,1){23}}
\put(25,50){\line(-1,-1){23}}
\put(50,25){\line(-1,-1){23}}
\put(25,50){\line(1,0){23}}
\put(50,25){\line(0,1){23}}
\put(0,0){\circle*{5}}
\put(50,25){\circle*{5}}
\put(25,50){\circle*{5}}
\put(0,25){\circle{5}}
\put(25,0){\circle{5}}
\put(50,50){\circle{5}}
\put(0,0){\line(-1,-1){10}}
\put(50,25){\line(1,0){10}}
\put(25,50){\line(0,1){10}}
\put(-2,25){\line(-1,0){10}}
\put(25,-2){\line(0,-1){10}}
\put(52,52){\line(1,1){10}}
\end{picture}
}}}
\thicklines
\multiput(0,-25)(-75,75){2}{
\multiput(50,25)(20,10){8}{\put(2,0){\line(2,1){10}}}}
\multiput(0,-25)(150,75){2}{
\multiput(50,25)(-21,21){4}{\put(2,0){\line(-1,1){11}}}}
\end{picture}
\caption{\label{fig:honeycomb graph}
\textit{A fundamental parallelogram of the lattice $3\gL$ and a piece of the bipartite graph $\gG$ for
$\gA=\{(1,0),\,(0,1),\,(-1,-1)\}\subset\ZZ^2$, and  all 
$c_{\mf{a}}=1$. By identifying opposite sides of the parallelogram one
obtains the graph $\gG_3$.}}
\end{figure}
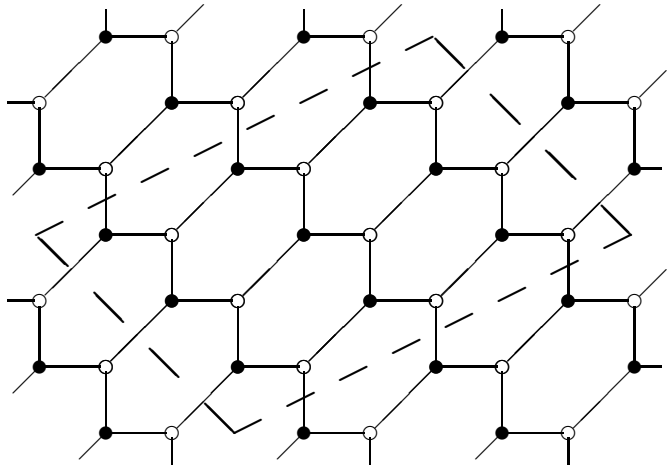

One can then construct a weighted bipartite graph $\gG$ as follows.
Bipartite graphs have two kinds of vertices,
often called black and white.
The set of \emph{black vertices} of $\gG$ is $\gL$. The set of \emph{white vertices} of $\gG$ is $\gA+\gL$. Note that $\gA+\gL$ is just one single coset of $\gL$ in $\ZZ^n$. In $\gG$ there is an (oriented) edge from vertex $\mf{v}_1$ to vertex $\mf{v}_2$ if and only if $\mf{v}_1$ is black,
$\mf{v}_2$ is white and $\mf{v}_2-\mf{v}_1\in\gA$. 
If $\mf{v}_2-\mf{v}_1=\mf{a}\in\gA$ the edge is said to be of \emph{type} $\mf{a}$ and gets \emph{weight} $c_\mf{a}$.

The graph $\gG$ is $\gL$-periodic and for every $N\in\NN$ one has the finite graph $\gG_N:=\modquot{\gG}{N\gL}$, which is naturally embedded in the torus $\modquot{\RR^n}{N\gL}$.
By a \emph{closed path} of length $2k$ on $\gG$ or $\gG_N$ we mean a sequence of edges 
$(e_1,e_2,\ldots,e_{2k-1},e_{2k})$ such that for $i=1,\ldots,k$
the intersection
$e_{2i-1}\cap e_{2i}$ contains a white vertex
and $e_{2i}\cap e_{2i+1}$ contains a black vertex; here $e_{2k+1}=e_1$. 
By the weight of such a path we mean the product of the weights of the edges $e_1,e_2,\ldots,e_{2k-1},e_{2k}$. 
We denote the set of closed paths of length $2k$ on $\gG_N$ by
$\gG_N(2k)$.
Enumerating the closed paths on $\gG_N$ according to length and weight
we prove in Section \ref{section:graph} that this leads to a new interpretation of the polynomials $B_N(z)$: 
\begin{theorem}\label{thm:BN graph}
\begin{equation}\label{eq:BN graph}
B_N(z)\,=\,z^{N^n}\exp\left(-\sum_{k\geq 1}\sum_{\gamma\in\gG_N(2k)}
\mathrm{weight}(\gamma)z^{-k}\right)\,.
\end{equation}
\qed
\end{theorem}
Formulas (\ref{eq:BN diff}) and (\ref{eq:BN graph}) transfer the enumerative data between the two models.

We pass to algebraic geometry with the Laurent polynomial
\begin{equation}\label{eq:W}
W(x_1,\ldots,x_n)\,=\,
\sum_{\mf{a},\mf{b}\in\gA}c_\mf{a}c_\mf{b}
\mf{x}^{\mf{a}-\mf{b}}\,\in\ZZ[x_1^{\pm 1},\dots,x_n^{\pm 1}]\,,
\end{equation}
which satisfies 
$|\FD(\mf{t})|^2 = W(e^{2\pi i t_1},\ldots,e^{2\pi i t_n})$; here
$\mf{x}^{\gl}:=\prod_{j=1}^n x_j^{\gl_j}$
if $\gl=(\gl_1,\ldots,\gl_n)\in\gL$.
For $N\in\NN$ let $\ru_N$ denote the group of $N$-th roots of unity and let 
$
\rul_N:=\mathrm{Hom}(\gL,\ru_N)
$
be the group of homomorphisms from the lattice $\gL$ to $\ru_N$. 
Thus the defining formula (\ref{eq:BN diff}) can be rewritten as:
\begin{equation}\label{eq:BN unity}
\sB_N(z)=\prod_{\mf{x}\in\rul_N}(z-W(\mf{x}))\,. 
\end{equation}
Written in the form (\ref{eq:BN unity}) the polynomials $\sB_N(z)$
appear as direct generalizations of quantities introduced by Lehmer 
\cite{L} for a $1$-variable (i.e. $n=1$) polynomial $W(x)$. 
Using (\ref{eq:BN unity}) one can easily show 
(Proposition \ref{prop:divisibility sequence}) that
the polynomials $\sB_N(z)$ have integer coefficients and that 
$\sB_{N'}(z)$ divides $\sB_N(z)$ in $\ZZ[z]$ if $N'$ divides $N$ 
in $\ZZ$. Thus 
for $h\in\ZZ$ also $\sB_N(h)$ is an integer.
Lehmer was particularly interested in the prime factorization of these integers in case $n=1$ \cite{L,EW,RV}. Also for general $n\geq 1$ these prime factorizations must be interesting, for instance because they relate to counting points on varieties over finite fields;
see Section \ref{section:geometry} for details.
\emph{Thus prime factorization gives a third occurrence of $\sB_N(z)$ in enumerative problems, related to counting points on varieties over finite fields.}
$\sQ(z)$ appears in \cite{Bo,De,RV} as \emph{Mahler measure} with ties to special values of $L$-functions. It would be nice if the limit formula (\ref{eq:radius})
together with the prime factorization of the numbers $\sB_N(z)$ (with $z\in\ZZ$) could shed new light on these very intriguing ties.

In \emph{Section \ref{section:diffraction}} we study the density distribution
of the level sets in the diffraction pattern. Passing from measures to complex functions with the Hilbert transform we find one interpretation of $\sQ(z)$, $\sB_N(z)$ and (\ref{eq:radius}).
In \emph{Section \ref{section:Laplace}} we briefly discuss another interpretation in connection with the spectrum of a discretized Laplace operator. In \emph{Section \ref{section:graph}} we prove 
Theorem \ref{thm:BN graph}.
In \emph{Section \ref{section:geometry}} we pass to toric geometry, where the diffraction pattern reappears as the intersection of a real torus with a family of hypersurfaces in a complex torus and where
$\log\sQ(z)$ becomes a period integral, while the prime factorization of $\sB_N(z)$ for $z\in\ZZ$ somehow relates to counting points on those hypersurfaces over finite fields.
In \emph{Section \ref{section:moments}} we discuss sequences of integers which appear as moments of measures, path counts on graphs and  coefficients in Taylor expansions of solutions of Picard-Fuchs differential equations.
Finally, in \emph{Section \ref{section:examples}} we present some concrete examples. 

\section{The diffraction pattern.}\label{section:diffraction}

The function $|\FD(\mf{t})|^2$ is periodic with period lattice $\vL$ dual to the lattice $\gL$:
$$
\begin{array}{lcl}
\vL&:=&\{ \mf{t}\in\RR^n\;|\; \langle \mf{t},\mf{a}-\mf{b}\rangle\in\ZZ\,,\;
\forall \mf{a},\mf{b}\in\gA\}\,,
\\
\gL&:=&\ZZ\mathrm{-Span}\{ \mf{a}-\mf{b}\;|\; \mf{a},\mf{b}\in\gA\}.
\end{array}
$$
\emph{Throughout this note we assume that the lattices $\gL$ and $\vL$ have rank $n$.}

Because of this periodicity $|\FD(\mf{t})|^2$ descends to a function on $\modquot{\RR^n}{\vL}$.
Defining for $\mf{t}\in\RR^n$ the function 
$e_{\mf{t}}:\RR^n\rightarrow\CC$ by
$e_{\mf{t}}(\mf{v})=e^{2\pi i\langle \mf{t},\mf{v}\rangle}$
we obtain an isomorphism of real tori
\begin{equation}\label{eq:tori}
\modquot{\RR^n}{\vL}\;\simeq\; \UL\,,\qquad \mf{t}\mapsto e_{\mf{t}}\,,
\end{equation}
where $\UL:=\mathrm{Hom}(\gL,\UU)$ is the torus of group homomorphisms from the lattice $\gL$ to the unit circle 
$\UU:=\{x\in \CC\:|\;|x|=1\,\}$.
Recall that a group homomorphism $\psi:\gL\rightarrow\UU$ induces an algebra homomorphism $\psi_*$ from the group algebra $\CC[\gL]$ to $\CC$.
Thus $\CC[\gL]$ is the natural algebra of functions on $\UL$ and
$\psi_*$ evaluates functions at the point $\psi$ of $\UL$.
The inclusion $\gL\subset\ZZ^n$ identifies $\CC[\gL]$
with the subalgebra of the algebra of Laurent polynomials
$\CC[x_1^{\pm 1},\dots,x_n^{\pm 1}]$, which consists of $\CC$-linear
combinations of the monomials $\mf{x}^{\gl}:=\prod_{j=1}^n x_j^{\gl_j}$
with $\gl=(\gl_1,\ldots,\gl_n)\in\gL$.
Thus,
via (\ref{eq:tori}), the function $|\FD(\mf{t})|^2$ coincides with  the Laurent polynomial $W(x_1,\ldots,x_n)$ defined in (\ref{eq:W}).

Positivity of the coefficients $c_\mf{a}$ implies that the function $|\FD(\mf{t})|^2$ attains its maximum exactly at the points
$\mf{t}\in\vL$. In terms of the torus $\UL$ and the function $W:\UL\rightarrow\RR$ this means that $W$ attains its maximum exactly at the origin $\mf{1}$ of the torus group $\UL$:
$$
\forall\mf{x}\in\UL\setminus\{\mf{1}\}:\quad W(\mf{x})\,<\,W(\mf{1})\,=\,C^2\qquad\textrm{with}\quad
C:=\sum_{\mf{a}\in\gA}c_\mf{a}\,.
$$

Some important aspects of the density distribution in the diffraction pattern are captured by the function
\begin{equation}
\label{eq:volume1}
V:\RR\rightarrow\RR\,,\qquad
V(r):=\frac{\mathrm{volume}
\{\mf{t}\in \modquot{\RR^n}{\vL}\;|\; |\FD(\mf{t})|^2\leq r\}}
{\mathrm{volume}(\,\modquot{\RR^n}{\vL}\,)}\,.
\end{equation}
We view the derivative $dV(r)$ of $V$ as a measure on $\RR$. 
In our analysis it will be important that
the measure $dV(r)$ is also
the push forward of the standard measure $dt_1\,dt_2\,\dots dt_n$ on 
$\RR^n$ by the function $|\FD (\mf{t})|^2$.

Another insight into the diffraction pattern comes from its intersection with the torsion subgroup of $\UL$.
For $N\in\NN$ let $\ru_N\subset\UU$ denote the group of $N$-th roots of unity. Then the group of $N$-torsion points in $\UL$ is
$\rul_N:=\mathrm{Hom}(\gL,\ru_N)$ and (\ref{eq:multiplicities})
can be rewritten as
$$
\mul_N(r)=\sharp (\rul_N\,\cap\,W^{-1}(r))
\qquad \textrm{for}\quad N\in\NN,\,r\in\RR\,.
$$
Moreover we set, in analogy with (\ref{eq:volume1}),
$$
V_N(r):=\frac{1}{N^n}\sharp
\{\mf{x}\in \rul_N\;|\; W(\mf{x})\leq r\}
\qquad \textrm{for}\quad N\in\NN,\,r\in\RR\,.
$$
The derivative of the step function $V_N:\RR\rightarrow\RR$ 
is the distribution
\begin{equation} \label{eq:dVN}
dV_N (r)\,=\,N^{-n}\sum_r \mul_N(r)\,\delta_r\,,
\end{equation}
which assigns to a continuous function $f$ on $\RR$ the value
\begin{equation} \label{eq:measure dVN}
\int_\RR f(r)dV_N (r):= N^{-n}\sum_r \mul_N(r)\,f(r)
= N^{-n}\sum_{\mf{x}\in \rul_N}
f(W(\mf{x}))\,.
\end{equation}
One thus finds a limit of distributions
\begin{equation}
\label{eq:measure limit}
\lim_{N\to\infty} dV_N(r)\,=\,dV(r)\,;
\end{equation}
by definition, this means that 
for  every continuous function $f$ on $\RR$
\begin{equation}
\label{eq:Riemann sum} 
\lim_{N\to\infty}\int_\RR f(r)dV_N (r)\,=\,\int_\RR f(r)dV (r)\,.
\end{equation}

Measure theory is connected with complex function theory by the Hilbert transform. The Hilbert transform of the measure $dV(r)$ is the function
$-\frac{1}{\pi}\sH (z)$ defined by
\begin{equation}
\label{eq:hilbert transform}
\sH (z)\,:=\,\int_\RR \frac{1}{z-r} dV(r)\qquad\textrm{for}\quad 
z\in\CC\setminus\cI\,;
\end{equation}
here $\cI:=\overline{\{r\in\RR\:|\:0<V(r)<1\}}\,\subset\,[0,C^2]$ is the support of the measure $dV(r)$.
The measure can be recovered from its Hilbert transform because
for $r_0\in\RR$
$$
\frac{dV}{dr}(r_0)=\frac{1}{2\pi i}\lim_{\ep\in\RR,\ep\downarrow 0}
\left(\sH (r_0+i\ep)\,-\,\sH (r_0-i\ep) \right)\,.
$$
Another way of writing the connection between $dV(r)$ and $\sH(z)$ is 
\begin{equation}
\label{eq:recover2}
\frac{1}{2\pi i}\oint_\gamma f(z)\sH (z)dz\,=\,
\int_\RR\left(\frac{1}{2\pi i}\oint_\gamma \frac{f(z)}{z-r} 
dz\right)dV(r)\,=\,
\int_\RR f(r)dV(r)
\end{equation}
for
holomorphic functions $f$ defined on some open neighborhood $U$
of the interval $\cI$ in $\CC$ and closed paths $\gamma$ in 
$U\setminus\cI$ encircling $\cI$ once counter clockwise.

Next we consider the function
\begin{equation}
\label{eq:Q}
\sQ:\CC\setminus\cI\longrightarrow\RR_{>0}\,,\qquad
\sQ(z):=\exp\left(-\int_\RR\log|z-r|dV(r)\right)\,.
\end{equation}
This function satisfies
$
\frac{d}{dz}\log\sQ(z)=-\sH(z)
$
and thus (\ref{eq:recover2}) can be rewritten as
$$
\frac{-1}{2\pi i}\oint_\gamma f(z)d\log\sQ(z)\,=\,
\int_\RR f(r)dV(r)\,.
$$
\emph{This means that, at least intuitively, the functions
$\sQ(z)$ and $e^{-V(r)}$ correspond to each other via some kind of comparison isomorphism.}

In order to find the analogue of (\ref{eq:measure limit}) in terms of functions on $\CC\setminus\cI$
we apply (\ref{eq:measure dVN})
to the function $f(r)=\log|z-r|$ on $\RR$ with fixed $z\in\CC\setminus\cI$:
\begin{equation} \label{eq:measure dVN2}
\int_\RR \log|z-r|dV_N (r)= N^{-n}\log|\sB_N(z)|
\end{equation}
where 
$ 
\sB_N(z)=\prod_{r\in\cI}(z-r)^{\mul_N(r)}=
\prod_{\mf{x}\in\rul_N}(z-W(\mf{x}))
$
as in (\ref{eq:BN diff}) and (\ref{eq:BN unity}).

Combining (\ref{eq:Riemann sum}), (\ref{eq:Q}) and (\ref{eq:measure dVN2}) we find the limit announced in (\ref{eq:radius}):

\begin{proposition}\label{prop:Mahler limit}
$\displaystyle{\quad
\sQ(z)=\lim_{N\to\infty} |\sB_{N}(z)|^{-N^{-n}}
}$
for every $z\in\CC\setminus\cI$.
\qed
\end{proposition}

\

\

\section{The Laplacian perspective.}
\label{section:Laplace}

Convolution with the distribution $\cD$ gives the operator
$\sD f(\mf{v}):=\displaystyle{
\sum_{\mf{a}\in\gA}c_\mf{a} f(\mf{v}-\mf{a})}$
on the space of $\CC$-valued functions on $\RR^n$. Let $\osD f(\mf{v}):=\displaystyle{
\sum_{\mf{a}\in\gA}c_\mf{a} f(\mf{v}+\mf{a})}$ and
$$
\gD:=\sD\osD\,,\qquad
\gD f (\mf{v})=
\sum_{\mf{a},\mf{b}\in\gA}c_\mf{a}c_\mf{b} f(\mf{v}+\mf{a}-\mf{b})\,.
$$

\

For a sufficiently differentiable function $f$ on $\RR^n$ the Taylor expansion 
$$
\gD f(\mf{v})
=C^2f(\mf{v})\,+\,\textstyle{\frac{1}{2}}\;\displaystyle{
\sum_{i,j=1}^n
\sum_{\mf{a},\mf{b}\in\gA}c_\mf{a}c_\mf{b} (a_i-b_i)(a_j-b_j)
\frac{\partial^2 f}{\partial v_i \partial v_j}(\mf{v})\,+}\ldots
$$
shows that the difference operator
$\gD-C^2$ is a \emph{discrete approximation of the Laplace operator} corresponding to the Hessian 
of the function $|\FD(\mf{t})|^2$ at its maximum.

\

\noindent
\textbf{Remark}.
In \cite{GKT} Gieseker, Kn\"orrer and Trubowitz investigate 
Schr\"odinger equations in solid state physics via a discrete approximation of the Laplacian. In their situation the Schr\"odinger operator is the discretized Laplacian \emph{plus a periodic potential function}. So from the perspective of \cite{GKT}
the present note deals only with the (simple) case of zero potential.
On the other hand we consider more general discretization schemes and possibly higher dimensions.

\

We now turn to the spectrum of $\gD$.
For $\mf{t}\in\RR^n$ the function 
$e_{\mf{t}}:\RR^n\rightarrow\CC$ given by
$e_{\mf{t}}(\mf{v})=e^{2\pi i\langle \mf{t},\mf{v}\rangle}$ is an eigenfunction for $\gD$ with eigenvalue $|\FD(\mf{t})|^2$:
$$
\gD e_{\mf{t}}(\mf{v})\, =\,
\sum_{\mf{a},\mf{b}\in\gA}c_\mf{a}c_\mf{b} 
e^{2\pi i\langle \mf{t},\mf{v}+\mf{a}-\mf{b}\rangle}\, =\,
|\FD(\mf{t})|^2\:e_{\mf{t}}(\mf{v})\,.
$$
Take a positive integer $N$.
The space of $\CC$-valued
$C^\infty$-functions on $\RR^n$ which are periodic for the sublattice $N\gL$ of $\gL$ is spanned by the functions
$e_{\mf{t}}$ with $\mf{t}$ in the dual lattice
$\frac{1}{N}\vL$.
The characteristic polynomial of the restriction of $\gD$ to this space is therefore (see (\ref{eq:BN diff}) and (\ref{eq:BN unity}))
$$
\prod_{\textstyle{\mf{t}}\in\modquot{\textstyle{\frac{1}{N}}\vL}{\vL}}
(z-|\FD(\mf{t})|^2)\;=\;
\prod_{\mf{x}\in\rul_N}(z-W(\mf{x}))\;=\;
\prod_{r\in\RR}(z-r)^{\mul_N(r)}\;=\;\sB_N(z)\,.
$$

\emph{With (\ref{eq:dVN}) and (\ref{eq:measure limit}) the measure $dV(r)$ can now be interpreted as the density of the eigenvalues of $\gD$ on the space of $\CC$-valued 
$C^\infty$-functions on $\RR^n$ which are periodic for some  sublattice $N\gL$ of $\gL$.}

\section{Enumeration of paths on a periodic weighted bipartite graph;
proof of Theorem \ref{thm:BN graph}.}
\label{section:graph}

In this section we prove Theorem \ref{thm:BN graph}. Recall from the 
Introduction just before Theorem \ref{thm:BN graph} the various ingredients: the finite set $\gA\subset\ZZ^n$, the weights
$c_\mf{a}$, the lattice $\gL$ and the graphs $\gG$ and $\gG_N$. Recall also the closed paths on $\gG_N$, their lengths and weights, and 
the set $\gG_N(2k)$ of closed paths of length $2k$ on $\gG_N$.

Consider a path $(e_1,e_2,\ldots,e_{2k-1},e_{2k})$ on $\gG$
with edge $e_i$ going from black to white if $i$ is odd,
respectively from white to black if $i$ is even. Let $\mf{s}$ denote the starting point of the path (i.e. the black vertex of edge $e_1$). Let for $j=1,\ldots,k$ edge $e_{2j-1}$ be of type $\mf{a}_j$ and edge $e_{2j}$ of type $\mf{b}_j$.
Then the end point of the path (i.e. the black vertex of $e_{2k}$) is
$\mf{s}+\sum_{j=1}^k (\mf{a}_j-\mf{b}_j)$. The weight of the path is
$\prod_{j=1}^k c_{\mf{a}_j}c_{\mf{b}_j}$. The path closes on $\gG_N$ if and only if
$
\sum_{j=1}^k (\mf{a}_j-\mf{b}_j)\,\in\,N\gL
$.

Next recall from (\ref{eq:W}) that
$
W(\mf{x})\,=\,
\sum_{\mf{a},\mf{b}\in\gA}c_\mf{a}c_\mf{b}
\mf{x}^{\mf{a}-\mf{b}}
$
and set
\begin{equation}\label{eq:N moments1}
\mo_k^{(N)}\,:=\, 
N^{-n}\sum_{\mf{x}\in\rul_N}W(\mf{x})^k\,.
\end{equation}
So $\mo_k^{(N)}$ is the sum of the coefficients of those monomials in
$W(\mf{x})^k$ with exponent in $N\gL$. In view of the above considerations $\mo_k^{(N)}$ is therefore equal to the sum of the weights of the paths on $\gG$ which start at $\mf{s}$, have length $2k$
and close in $\gG_N$.
Since on $\gG_N$ there are $N^n$ black vertices and on a path of length $2k$ there are $k$ black vertices we conclude 
\begin{equation}\label{eq:mk graph}
\frac{N^n}{k}\mo_k^{(N)}\,=\,
\sum_{\gamma\in\gG_N(2k)}\mathrm{weight}(\gamma)\,.
\end{equation}
From (\ref{eq:BN unity}) one sees for $|z|>C^2$
\begin{equation}\label{eq:logBN} 
N^{-n}\log B_N(z)\,=\,\log z\,+
\,N^{-n}\sum_{\mf{x}\in\rul_N}\log(1-W(\mf{x})z^{-1})\,=\,
\log z-\sum_{k\geq 1}\frac{\mo_k^{(N)}}{k} z^{-k}\,.
\end{equation}
Combining (\ref{eq:mk graph}) and (\ref{eq:logBN}) we find
$$
B_N(z)\,=\,z^{N^n}\exp\left(-\sum_{k\geq 1}
\sum_{\gamma\in\gG_N(2k)}\mathrm{weight}(\gamma)\;
z^{-k}\right)\,.
$$
This finishes the proof of Theorem \ref{thm:BN graph}.\qed

\

\noindent
\textbf{Remark.}
In the Laplacian perspective $N^n\mo_k^{(N)}$ is the trace of the operator
$\gD ^k$ on the space of $\CC$-valued 
$C^\infty$-functions on $\RR^n$ which are periodic for the  sublattice $N\gL$ of $\gL$. The polynomial $B_N(z)$ is the characteristic polynomial
of $\gD$ on this space. Formula (\ref{eq:logBN}) gives the well-known relation between the characteristic polynomial of an operator and the traces of its powers.

\

\noindent
\textbf{Remark.}
One may refine the above enumerations by keeping track of the homology class to which the closed path belongs. That means that instead of (\ref{eq:N moments1}) one extracts from the polynomial $W(\mf{x})^k$
the subpolynomial consisting of terms with exponent in $N\gL$.
Such a refinement of the enumerations with homology data appears also in the theory of dimer models (cf. \cite{KOS}), but its meaning for the diffraction pattern is not clear.

\section{Algebraic geometry.}\label{section:geometry}

The polynomial $\sB_N(z)=\prod_{\mf{x}\in\rul_N}(z-W(\mf{x}))$ 
has coefficients in the ring of integers of the cyclotomic field
$\QQ(\ru_N)$ and is clearly invariant under the Galois group of
$\QQ(\ru_N)$ over $\QQ$. Consequently, the coefficients of $\sB_N(z)$
lie in $\ZZ$. The same argument applies to the polynomial
$\sB_N(z)\sB_{N'}(z)^{-1}=\prod_{\mf{x}\in\rul_N\setminus\rul_{N'}}
(z-W(\mf{x}))$ if $N'$ divides $N$.
Thus we have proved

\begin{proposition}\label{prop:divisibility sequence}
For every $N\in \NN$
the coefficients of $\sB_N(z)$
lie in $\ZZ$.
If $N'$ divides $N$ in $\ZZ$, then $\sB_{N'} (z)$ divides
$\sB_N (z)$ in $\ZZ[z]$.
\qed
\end{proposition} 

\

Fix a prime number $p$ and a positive integer
$\nu\in\ZZ_{>0}$.
Let $\WW(\FF_{p^\nu})$ denote the ring of Witt vectors of the finite field
$\FF_{p^\nu}$ (see e.g. \cite{B}). So, $\WW(\FF_{p^\nu})$ is a complete discrete valuation 
ring with
maximal ideal $p\WW(\FF_{p^\nu})$ and residue field $\FF_{p^\nu}$.
The Teichm\"{u}ller lifting is a map $\tau:\FF_{p^\nu}\longrightarrow\WW(\FF_{p^\nu})$ such that
$$
\quad x\equiv\tau (x)\bmod p\,,\qquad
\tau(xy)=\tau(x)\tau(y)\qquad\forall x,y\in\FF_{p^\nu}\,.
$$
Every non-zero $x\in\FF_{p^\nu}$ satisfies
$$
x^{p^\nu-1}=1\,.
$$
Thus there is an isomorphism $\ru_{p^\nu-1}\simeq \FF_{p^\nu}^*$.
Such an isomorphism composed with the Teichm\"{u}ller lifting gives an embedding $j:\ru_{p^\nu-1}\hookrightarrow\WW(\FF_{p^\nu})$.
Thus for $\mf{x}\in\rul_{p^\nu-1}$ we get 
$$
W(j(\mf{x}))\in\WW(\FF_{p^\nu})\,.
$$
Recall the $p$-adic valuation on $\ZZ$: for $k\in\ZZ$, $k\neq 0$:
$$
v_p(k)\,:=\,\max\{v\in\ZZ\;|\;p^v \textrm{ divides } k\}\,.
$$

\begin{proposition}
For $p,\nu$ as above and for $z\in\ZZ$ the $p$-adic valuation of the integer $\sB_{p^\nu-1} (z)$ satisfies
\begin{equation}\label{eq:pval BN}
v_p(\sB_{p^\nu-1} (z))\geq\sharp\{\xi\in (\FF_{p^\nu}^*)^n\;|\;
W(\xi)=z\; \textrm{in}\; \FF_{p^\nu}\;\}\,.
\end{equation}
\end{proposition}
\proof
From (\ref{eq:BN unity}) we obtain the product decomposition,
with factors in $\WW(\FF_{p^\nu})$,
$$
\sB_N(z)=\prod_{\xi\in (\FF_{p^\nu}^*)^n}(z-W(\tau(\xi)))\,. 
$$
The result of the proposition now follows because
$$
z-W(\tau(\xi))\in p\WW(\FF_{p^\nu})\qquad\Leftrightarrow\qquad
W(\xi)=z\; \textrm{in}\; \FF_{p^\nu}
$$
\qed

\

\noindent
\textbf{Remark.} In (\ref{eq:warning}) we give an example showing that 
in (\ref{eq:pval BN}) we may have a strict inequality. 

\

\noindent
\textbf{Remark about the relation with Mahler measure and L-functions.}\\
The \emph{logarithmic Mahler measure} $\mf{m} (F)$ and the \emph{Mahler measure} $\mf{M} (F)$ of a Laurent polynomial $F(x_1,\ldots,x_n)$ with complex coefficients are:
\begin{eqnarray*}
\mf{m} (F)&:=&\frac{1}{(2\pi i)^n}\oint\!\!\oint_{|x_1|=\ldots=|x_n|=1}
\log |F(x_1,\ldots,x_n)|\,
\frac{dx_1}{x_1}\cdot\ldots\cdot\frac{dx_n}{x_n}\,,\\
\mf{M} (F)&:=&\exp (\mf{m} (F))\;.
\end{eqnarray*}
Boyd \cite{Bo} gives a survey of many (two-variable) Laurent polynomials for which 
$\mf{m} (F)$ equals 
(numerically to many decimal places) a `simple' non-zero rational number times the derivative at $0$ of the L-function of the projective plane curve  $Z_F$ defined by the vanishing of $F$:
\begin{equation}\label{eq:special L}
\mf{m}(F)\cdot\QQ^*\;=\;\mathrm{L}'(Z_F,0)\cdot\QQ^*\,.
\end{equation}
Deninger \cite{De} and Rodriguez Villegas \cite{RV} showed that the experimentally observed relations (\ref{eq:special L}) agree with predictions from the Bloch-Beilinson conjectures.
Rodriguez Villegas \cite{RV} provided actual proofs for a few special examples.

Since the measure $dV(r)$ is 
the push forward of the measure $dt_1\,dt_2\,\dots dt_n$ on 
$\RR^n$ by the function $|\FD(\mf{t})|^2$, one can rewrite Formula
(\ref{eq:Q}) as:
\begin{equation}\label{eq:logQ}
-\log\sQ (z)\,=\,\frac{1}{(2\pi i)^n}\int_{\UU^n}
\log |z-W(x_1,\ldots,x_n)|\, \frac{dx_1}{x_1}\,\frac{dx_2}{x_2}\,\dots \,\frac{dx_n}{x_n}\,.
\end{equation}
On the right hand side of (\ref{eq:logQ}) we now recognize the
logarithmic Mahler measure of the Laurent polynomial 
\mbox{$z-W(x_1,\ldots,x_n)\,\in\,\CC[x_1^{\pm 1},\ldots,x_n^{\pm 1}]$.}

For fixed $z\in\ZZ$ Formulas (\ref{eq:radius}) and (\ref{eq:pval BN}) provide a link between $\sQ (z)$ and counting points over finite fields 
on the variety with equation $W(x_1,\ldots,x_n)=z$. It may be an interesting challenge to further extend these ideas to a proof of a result like (\ref{eq:special L}).

\section{Moments.}\label{section:moments}
Important invariants of the measure $dV(r)$ are its
\emph{moments} $\mo_k$ ($k\in\ZZ_{\geq 0}$):
\begin{equation}\label{eq:moments}
\begin{array}{lrl}
\mo_k&:=&\int_\RR r^k dV(r)\;=\;
 \int_0^1\ldots\int_0^1
|\FD(t_1,\ldots,t_n)|^{2k}dt_1\ldots dt_n \\
&=&\textrm{constant term of Fourier series }\; |\FD(t_1,\ldots,t_n)|^{2k}\\
&=&\textrm{constant term of Laurent polynomial }\; W(x_1,\ldots,x_n)^k\,.
\end{array}
\end{equation}
The relation between the moments and the functions $\sH(z)$, $\sQ(z)$ defined in (\ref{eq:hilbert transform}) and (\ref{eq:Q}) is:
for $z\in\RR,\,z>C^2$,
\begin{equation}
\label{eq:HQmoments}
\sH (z)\,=\,\sum_{k\geq 0} \mo_k z^{-k-1}
\,,\qquad
\sQ(z)=z^{-1}\exp\left(\sum_{k\geq 1}\frac{\mo_k}{k} z^{-k}\right)\,.
\end{equation}

It is clear that the moments $\mo_k$ of $dV(r)$ are non-negative integers. They
satisfy all kinds of arithmetical relations. There are, for instance, recurrences like (\ref{eq:recurrence}) and congruences like the following 

\begin{lemma}\label{lemma congruenties}
$
\mo_{kp^{\ga+1}}\equiv \mo_{kp^\ga}\bmod p^{\ga+1}
$
for every prime number $p$ and $k,\ga\in\ZZ_{\geq 0}$.
\end{lemma}
\proof
The Laurent polynomial $W(x_1,\ldots,x_n)$ has coefficients in $\ZZ$.
Therefore
$$
W(x_1,\ldots,x_n)^{kp^{\ga+1}}\equiv W(x_1^p,\ldots,x_n^p)^{kp^\ga}
\bmod p^{\ga+1}\ZZ[x_1^{\pm 1},\ldots,x_n^{\pm 1}]\,.
$$
The lemma follows by taking constant terms.
\qed

\

Theorems 1.1, 1.2, 1.3 in \cite{Ari} together with the above lemma immediately yield the following integrality result for series and product expansions:

\begin{corollary}\label{cor:Qseries}
For $z\in\RR,\,z>C^2$
\begin{equation}\label{eq:Qseries}
z^{-1}\exp\left(\sum_{k\geq 1}\frac{\mo_k}{k} z^{-k}\right)
\;=\;z^{-1}+\sum_{k\geq 1} A_k z^{-k-1}
\;=\; z^{-1}\prod_{k\geq 1}(1-z^{-k})^{-b_k}
\end{equation}
with $A_k,b_k\in\ZZ$ for all $k\geq 1$.
\qed
\end{corollary}

\

\noindent
\textbf{Remark}.
In \cite{Ari, EPW} the result of Corollary \ref{cor:Qseries} is used to interpret $z\sQ (z)$ as the Artin-Mazur zeta function of a dynamical system, provided the integers $b_k$ are not negative. We have not yet found such a dynamical system within the present framework.

\

For $N\in\NN$ the moments of the measure  $dV_N(r)$ are, by definition,
$$
\mo_k^{(N)}:=\int_\RR r^k dV_N(r)= 
N^{-n}\sum_{\mf{x}\in\ru_N}W(\mf{x})^k\,.
$$
These are the same numbers as in (\ref{eq:N moments1}).

\begin{proposition}\label{prop:moment limit}
With the above notations we have 
$$
\begin{array}{rcl}
\mo_k^{(N)}&\geq&\mo_k\;\geq\;0\qquad \textrm{for all}\qquad N,\,k\,,\\
\mo_k^{(N)}&=&\mo_k\qquad \textrm{if}\qquad 
N>k\max_{\mf{a},\mf{b}\in\gA}\max_{1\leq j\leq n}|a_j-b_j|\,.
\end{array}
$$
\end{proposition}
\proof
$N^n\,(\mo_k^{(N)}-\mo_k)$ is the sum of the coefficients of all non-constant monomials
in the Laurent polynomial $W(x_1,\ldots,x_n)^k$ with exponents divisible by $N$. Since all coefficients of $W(x_1,\ldots,x_n)$ are positive, this shows $\mo_k^{(N)}\geq \mo_k\geq 0$.

Assume $N>k\max_{\mf{a},\mf{b}\in\gA,\:1\leq j\leq n}|a_j-b_j|$.
Then all exponents in the monomials of the Laurent polynomial $W(x_1,\ldots,x_n)^k$ are $>-N$ and $<N$.
So only the exponent of the constant term is divisible by $N$.
Therefore
$\mo_k^{(N)}\,=\, \mo_k$.
\qed

\

Note the natural interpretation (and proof) of this proposition in terms of closed paths on the graph $\gG_N$: closed paths on $\gG_N$ which are too short are in fact projections of closed paths on $\gG$.

\

\begin{corollary}
For $\displaystyle{N>\ell\max_{\mf{a},\mf{b}\in\gA}
\max_{1\leq j\leq n}|a_j-b_j|}$ and $|z|>C^2$:
$$
\sQ(z)\cdot |B_N(z)|^{N^{-n}}
=\left|\exp\left(\sum_{k> \ell}\frac{\mo_k-\mo_k^{(N)}}{k} z^{-k}\right)\right|\,.
$$
This not only gives an estimate for the rate of convergence of
(\ref{eq:radius}) with respect to the usual absolute value on $\CC$, but it also yields the following congruence of power series in $z^{-1}$:
$$
B_N(z)^{-N^{-n}}\,\equiv\,
z^{-1}+\sum_{k\geq 1} A_k z^{-k-1}\;\bmod z^{-\ell-1}
$$
with $A_k$ as in (\ref{eq:Qseries})
\qed
\end{corollary}

\

\noindent
\textbf{Remark about the relation with the large complex structure limit.}
Since the measure $dV(r)$ is 
the push forward of the measure $dt_1\,dt_2\,\dots dt_n$ on 
$\RR^n$ by the function $|\FD(\mf{t})|^2$, one can rewrite  
(\ref{eq:hilbert transform}) as
$$
\sH (z)\,=\,
\frac{1}{(2\pi i)^n}\int_{\UU^n}
\frac{1}{z-W(x_1,\ldots,x_n)}\, \frac{dx_1}{x_1}\,\frac{dx_2}{x_2}\,\dots \,\frac{dx_n}{x_n}
$$
for $z\in\CC\setminus\cI$. From this (and the residue theorem) one sees that $\sH (z)$ is \emph{a period of some differential form of degree $n-1$
along some $(n-1)$-cycle on the hypersurface in $(\CC^*)^n$ given by the equation $W(x_1,\ldots,x_n)=z$.} As $z$ varies we get a $1$-parameter family of hypersurfaces. The function $\sH (z)$ is a solution of the Picard-Fuchs differential equation associated with (that 
$(n-1)$-form on) this family of hypersurfaces.
The Picard-Fuchs equation is equivalent with a recurrence relation for the coefficients $\mo_k$ in the power series expansion
(\ref{eq:HQmoments}) of $\sH (z)$ near $z=\infty$.
All this is standard knowledge about Calabi-Yau varieties near the large
complex structure limit and there is an equally standard algorithm to derive from the Picard-Fuchs differential equation enumerative information about numbers of instantons (or rational curves); see for instance \cite{Yau, M, S1}.

On the other hand, we have the enumerative data $\mo_k^{(N)}$ of the present paper. In the limit for $N\to\infty$ these yield the moments $\mo_k$ and, hence,
the Picard-Fuchs differential equation and eventually the instanton numbers.

\section{Examples}
\label{section:examples}
\subsection{$n=1$}\label{subsection:n=1}
Mahler measures of one variable polynomials have a long history with many interesting results; see the introductory sections of
\cite{Bo,EW,RV}.
We limit our discussion to one example, without a claim of new results. 
This simple, yet non-trivial, example has $n=1$, $\gA=\{-1,1\}\subset\ZZ$, $c_{-1}=c_1=1$ and hence
$$
|\FD(t)|^2=2+2\cos(4\pi t)\,,\qquad
W(x)=(x+x^{-1})^2\,.
$$
The moments are
$$
\mo_k\,=\,\textrm{constant term of } (x+x^{-1})^{2k}\,=\,
\left(\begin{array}{c} 2k\\ k\end{array}\right)
$$
and hence by (\ref{eq:HQmoments}): for $z\in\RR,\,z>4$,
\begin{eqnarray*}
\sH (z)&=&
\sum_{k\geq 0}\left(\begin{array}{c} 2k\\ k\end{array}\right)\;z^{-k-1}
\quad =\quad
\frac{1}{\sqrt{z(z-4)}}\,,\\
\sQ (z)&=&\exp\left(-\int\frac{dz}{\sqrt{z(z-4)}}\right)\,=\,
\frac{1}{2}\left(\,z-2-\sqrt{z(z-4)}\,\right)\,.
\end{eqnarray*}
Applying Formula (\ref{eq:N moments1}) to the present example we find
$$
\mo_k^{(N)}=\sum_{j\equiv k\bmod N}
\left(\begin{array}{c} 2k\\ j\end{array}\right)\,,
$$
which nicely illustrates Proposition \ref{prop:moment limit}.

Setting $z=2+u+u^{-1}$ one finds for the polynomials $\sB_N(z)$ defined in (\ref{eq:BN diff}):
\begin{eqnarray*} 
\sB_N(z)&=&\displaystyle{\prod_{x^2\in\ru_N}(z-2-x^2-x^{-2})\;=\;u^{-N}\prod_{x^2\in\ru_N}(u-x^2)(u-x^{-2})}\\[1ex]
\nonumber
&=&u^N+u^{-N}-2\\
&=&\left(\frac{1}{2}\left(\,z-2-\sqrt{z(z-4)}\,\right)\right)^N\,+\,
\left(\frac{1}{2}\left(\,z-2+\sqrt{z(z-4)}\,\right)\right)^N\,-\,2\\[1ex]
&=& -2+
2^{1-N}\sum_{j}\left(\begin{array}{c}N\\ 2j\end{array}\right)
z^j(z-4)^j(z-2)^{N-2j}\,.
\end{eqnarray*}
\emph{So, $\sB_N(z)$ is up to some shift and normalization the $N$-th \v{C}eby\v{s}ev polynomial.}

The above computation also shows
$
\sB_N(z)\,=\,\sQ (z)^N\,+\,\sQ (z)^{-N}-2
$
and thus, in agreement with (\ref{eq:radius}),
$$
\lim_{N\to\infty}\sB_N(z)^{-N^{-1}}\,=\,\sQ (z)\,.
$$

For actual computation of $\sB_N(z)$ in case $z\in\ZZ$ one can use the generating series identity:
$$
\sum_{N\geq 1} \sB_N(z)\,\frac{T^N}{N}\,=\,-\log\left(1-(z-4)\frac{T}{(1-T)^2}\right)\,.
$$
For $z=6$ one finds (using PARI)
\begin{eqnarray*}
\sum B_N(6)T^N&=&
2T + 12T^2 + 50T^3 + 192T^4 + 722T^5 + 2700T^6 + 10082T^7 \\
&& + 37632T^8 + 140450T^9 + 524172T^{10} + 1956242T^{11}  \\
&& + 7300800T^{12} + 27246962T^{13} + 101687052T^{14}  \\
&& + 379501250T^{15} + 1416317952T^{16} + 5285770562T^{17} + \ldots
\end{eqnarray*}
For primes $p$ in the displayed range the number $B_{p-1}(6)$ is
divisible by $p^2$ for $p\equiv\pm 1\bmod 12$ and is not
divisible by $p$ for $p\equiv\pm 5\bmod 12$ and is exactly divisible
by $p$ if $p=2,3$. 
We also checked $5^2\,|\,B_{24}(6)$ and $7^2\,|\,B_{48}(6)$.
This agrees with the number of solutions of the equation
$u+u^{-1}=4$ in $\FF_p$ and $\FF_{p^2}$.

If $z\in\ZZ, z>4$, then $\sQ (z)\,=\,\frac{1}{2}(\,z-2-\sqrt{z(z-4)})$ is a unit in the real quadratic field $\QQ(\sqrt{z(z-4)})$. According to Dirichlet's class number formula it relates to the $L$-function of this real quadratic field:
$$
\log(\sQ (z))\,=\,\frac{\sqrt{D}}{2h}L(1,\chi)
$$
where $D,\,h,\,\chi$ are the discriminant, class number, character, respectively,
of the real quadratic field $\QQ(\sqrt{z(z-4)})$ (see e.g. \cite{BS}).
The relations between Mahler measures and values of $L$-functions,
which have been observed for some curves, are perfect analogues 
of the above class number formula (see \cite{RV}).

\subsection{The honeycomb pattern.}\label{subsection:hexagon}
For a nice two-dimensional example we take
$\gA=\{(1,0),\,(0,1),\,(-1,-1)\}\subset\ZZ^2$, $c_{(1,0)}=c_{(0,1)}=c_{(-1,-1)}=1$ and hence
\begin{eqnarray*}
|\FD(t_1,t_2)|^2&\hspace{-.6em}=&\hspace{-.6em} 
3+2\cos(2\pi (t_1-t_2))+2\cos(2\pi (2t_1+t_2))+
2\cos(2\pi (t_1+2t_2)),\\
W(x_1,x_2)&\hspace{-.6em}=&\hspace{-.6em}
(x_1+x_2+x_1^{-1}x_2^{-1})(x_1^{-1}+x_2^{-1}+x_1x_2)
\\
&\hspace{-.6em}=&\hspace{-.6em}
x_1x_2^{-1}+x_1^2x_2+x_1^{-1}x_2+x_1x_2^2+
x_1^{-2}x_2^{-1}+x_1^{-1}x_2^{-2}+3\,.
\end{eqnarray*}
As a basis for the lattice $\gL$ we take $(2,1)$ and $(-1,-2)$. This leads to coordinates $u_1=x_1^2x_2$ and $u_2=x_1^{-1}x_2^{-2}$ on the torus $\UL$.
In these coordinates the function $W$ reads
\begin{eqnarray}\label{eq:n=2 Wu}
W(u_1,u_2)&=&u_1+u_1^{-1}+u_2+u_2^{-1}+u_1^{-1}u_2+u_1u_2^{-1}+3\\
\nonumber
&=&(u_1+u_2+1)(u_1^{-1}+u_2^{-1}+1)\,.
\end{eqnarray}
Figure \ref{fig:honeycomb graph} shows a piece of the graph $\gG$.
Figure \ref{fig:diffraction pattern} shows some level sets of the function $|\FD(t_1,t_2)|^2$.
The dual lattice $\vL$ is spanned by $(1,0)$ and 
$(\frac{1}{3},\frac{1}{3})$. The maximum of the function $|\FD(t_1,t_2)|^2$ equals $9$ and is attained at the points of $\vL$.
The minimum of the function $|\FD(t_1,t_2)|^2$ equals $0$ and is attained at the points of $(0,-\frac{1}{3})+\vL$ and $(\frac{1}{3},0)+\vL$.
There are saddle points with critical value $1$ at $(-\frac{1}{6},\frac{1}{3})+\vL$,
$(\frac{1}{3},-\frac{1}{6})+\vL$ and $(\frac{1}{6},\frac{1}{6})+\vL$. In terms of the coordinates $u_1,\,u_2$ the maximum lies at $(u_1,\,u_2)=(1,1)$, the minima at 
$(e^{2\pi i/3},e^{4\pi i/3}),\,(e^{4\pi i/3},e^{2\pi i/3})$    and the saddle points at
$(1,-1),(-1,1),(-1,-1)$.
The algebraic geometry of this example concerns the $1$-parameter family of elliptic curves with equation $z-W(u_1,u_2)=0$. In homogeneous coordinates $(U_0:U_1:U_2)$ on the projective plane $\PP^2$, with 
$u_1=U_1U_0^{-1}$, $u_2=U_2U_0^{-1}$, this becomes a homogeneous equation of degree $3$:
\begin{equation}\label{eq:elliptic pencil}
(U_0U_1+U_0U_2+U_1U_2)(U_0+U_1+U_2)-zU_0U_1U_2\,=\,0\,.
\end{equation}

Beauville \cite{Be} showed that there are exactly six semi-stable families of elliptic curves over $\PP^1$ with four singular fibres.
The pencil (\ref{eq:elliptic pencil}) is one of these six.
It has singular fibres at $z=0,1,9,\infty$ with Kodaira types
$I_2,I_3,I_1,I_6$, respectively. Note that the first three match the critical points and levels in the diffraction pattern.
After blowing up the points $(1,0,0),(0,1,0),(0,0,1)$ of $\PP^2$
one gets the DelPezzo surface $\textrm{dP}_3$.
The elliptic pencil (\ref{eq:elliptic pencil}) naturally lives on
$\textrm{dP}_3$. It has six base points, corresponding to six sections of the pencil. Since the base points have a zero coordinate, these sections do not intersect the real torus $\UL$.
Equations  (\ref{eq:n=2 Wu}) and (\ref{eq:elliptic pencil}) also appear in the literature in connection with the string theory of $\textrm{dP}_3$.

Formula (\ref{eq:moments}) and some manipulations of binomials give the moments:
$$
\mo_k\,=\,\sum_{j=0}^k 
\left(\begin{array}{c}k\\ j\end{array}\right)^2
\left(\begin{array}{c}2j\\ j\end{array}\right)\,.
$$
These numbers satisfy the recurrence relation (see \cite{SB} Table 7)
\begin{equation}\label{eq:recurrence}
(k+1)^2\mo_{k+1}=(10k^2+10k+3)\mo_k-9k^2\mo_{k-1}\,.
\end{equation}
We refer to \cite{SB} Example $\mathcal{C}$ and to \cite{S1} Example $\sharp 6$ for relations of these numbers to modular forms and instanton counts. Golyshev \cite{Go} can derive the recurrence
(\ref{eq:recurrence}) from the quantum cohomology of $\textrm{dP}_3$.
The numerical evidence for the relation 
(\ref{eq:special L}) between Mahler measure and L-function in this example is given in \cite{Bo} Table 2.

With Formulas (\ref{eq:N moments1}) and (\ref{eq:n=2 Wu}) one easily calculates
$$
\mo_k^{(N)}=\sum_{i_1\equiv i_2\bmod N,\,j_1\equiv j_2\bmod N}
\left(\begin{array}{c} k\\ i_1\end{array}\right)
\left(\begin{array}{c} k-i_1\\ j_1\end{array}\right)
\left(\begin{array}{c} k\\ i_2\end{array}\right)
\left(\begin{array}{c} k-i_2\\ j_2\end{array}\right)
\,.
$$
Note that these formulas confirm $\mo_k^{(N)}=\mo_k$ for $k<N$. 

Equation (\ref{eq:elliptic pencil}) is clearly invariant under 
permutations of $U_0,U_1,U_2$. Therefore the diffraction pattern has this $S_3$-symmetry too. Since only the critical points have a non-trivial stabilizer in $S_3$ the multiplicities $\mul_N(r)$ in this example satisfy
\begin{equation}\label{eq:n=2 multiplicities}
\begin{array}{lcll}
\mul_N(9)&=& 1& \forall N \\
\mul_N(0)&=& 2& \textrm{if}\quad 3|N \\
\mul_N(1)&\equiv& 3\bmod 6\quad & \textrm{if}\quad 2|N \\
\mul_N(r)&\equiv& 0\bmod 6\quad & \textrm{if}\quad r\neq 0,\,1,\,9,
\quad \forall N. 
\end{array}
\end{equation}
We have computed the numbers $\mul_N(r)$ for some values of $N$.
We found for instance
$$
\sB_{6}(z)=z^2\,(z-1)^{15}\,(z-3)^6\,(z-4)^6\,(z-7)^6\,(z-9)\,.
$$
We computed 
$W(u_1,u_2)\,=\,(u_1+u_2+1)(u_1^{-1}+u_2^{-1}+1)$ for $u_1,u_2\in\FF_7^*$:
the $(i,j)$-entry of the following $6\times 6$-matrix is $W(i,j)\bmod 7$:
$$
\left[
\begin{array}{rrrrrr}
2&3&0&3&0&1\\
3&3&4&0&1&1\\
0&4&0&1&4&1\\
3&0&1&3&4&1\\
0&1&4&4&0&1\\
1&1&1&1&1&1
\end{array}
\right]\,.
$$
This yields the following count of points over $\FF_7$:
$$
\begin{array}{rcrrrrrrrr}
z\bmod 7&:&0&1&2&3&4&5&6\\
\sharp\{\xi\in (\FF_7^*)^2\;|\;W(\xi)=z\; \textrm{in}\; \FF_7\;\}
&:&8&15&1&6&6&0&0
\end{array}
$$
Thus we see that the inequality in (\ref{eq:pval BN}) can be strict:
\begin{equation}\label{eq:warning}
v_7(\sB_6 (53))\,=\,12\,>\,6\,=\,\sharp\{\xi\in (\FF_7^*)^2\;|\;
W(\xi)=53\; \textrm{in}\; \FF_7\;\}\,.
\end{equation}

\

\

\textbf{Acknowledgement.} It is my pleasure to dedicate this paper to
Jaap Murre and Spencer Bloch on the occasion of their $75^{\textrm{th}}$,
respectively, $60^{\textrm{th}}$ birthdays. Both have been very important for my formation as a mathematician, from PhD-student time till present.

\end{document}